\documentclass[journal]{IEEEtran}

\usepackage{cite}

\ifCLASSINFOpdf
   \usepackage[pdftex]{graphicx}
\else
\fi

\usepackage{amsmath}

\begin{document}

\title{Distribution Network Marginal Costs --- Part II: Case Study Based Numerical Findings}

\author{Panagiotis~Andrianesis,
        and~Michael~Caramanis,~\IEEEmembership{Senior~Member,~IEEE}
\thanks{P. Andrianesis and M. Caramanis are with the Division of Systems Engineering, Boston University, 
Boston, MA, 02446 USA, e-mails: panosa@bu.edu, mcaraman@bu.edu. Research partially supported by the Sloan Foundation under grant G-2017-9723 and NSF AitF grant 1733827.}
}

\maketitle

\begin{abstract}
This two-part paper considers the day-ahead operational planning problem of a radial distribution network hosting Distributed Energy Resources (DERs), such as Solar Photovoltaic (PV) and Electric Vehicles (EVs).
In Part I, we develop a novel AC Optimal Power Flow (OPF) model that estimates dynamic Distribution nodal Location Marginal Costs (DLMCs) of real and reactive power including transformer degradation.
These all-inclusive DLMCs are key in the co-optimization of the distribution network and DER schedules.
Part II discusses the implication of using DLMCs to represent the benefit of shifting real/reactive power across time and achieve optimal Distribution Network and DER operation.
Moreover, it evaluates, elaborates and analyzes the novel AC OPF through extensive numerical results of actual distribution feeder-based case studies involving a wealth of future EV and PV adoption scenarios.
Optimal schedules are compared to schedules that rely on apparently reasonable approaches, namely Business as Usual (i.e., do nothing), Time of Use, and traditional line loss minimization. 
The overwhelming evidence of extensive numerical results supports the significant benefits of internalizing short-run marginal asset --- primarily transformer --- degradation.
Apart from optimal short-run scheduling, the proposed approach can harvest otherwise idle DER (PV/EV) reactive power compensation capabilities, increase distribution network DER (PV/EV) hosting capacity, and mitigate investments in distribution infrastructure that would be otherwise required to support distribution utilities' obligation to serve.
Finally, numerical evidence on the benefits of AC OPF modeling intertemporal transformer life degradation suggests that it may be worth considering the introduction of intra-day markets.
\end{abstract}


\IEEEpeerreviewmaketitle

\section{Introduction}

\IEEEPARstart{U}{nlike} a considerable body of the literature focusing on optimal open-loop DER scheduling, we focus on the derivation of comprehensive and accurate Distribution network Locational Marginal Costs (DLMCs) of real and reactive power, and employ them to co-optimize network decision trajectories along with the schedules of Distributed Energy Resources (DERs).
As such, we derive a mutually optimal coordination of network and DER decision trajectories that are fully adapted to each other through the equilibrium DLMCs.

Part I \cite{PartI} of this two-part paper focused on an innovative AC OPF framework that models network costs in a comprehensive manner, including the intertemporally coupled degradation of transformer life.
Furthermore, Part I focused on understanding and evaluating the components and building blocks that constitute DLMCs.
The contributions of Part II are as follows.

First, it relates DLMCs to financial incentives representing the benefit of shifting real/reactive power across time, and demonstrates how DLMCs convey sufficient information for achieving optimal Distribution Network and DER operation.

Second, it evaluates, elaborates and analyzes the novel AC OPF of Part I through extensive numerical findings of actual distribution feeder-based case studies involving a wealth of future EV and PV adoption scenarios.
The optimal schedules are compared to schedules that rely on conventional approaches that are often considered reasonable, including Business as Usual, 
time varying LMPs that do not vary across distribution nodes, and traditional Line Loss minimization.
Overwhelming evidence of extensive numerical results supports the significant benefits of internalizing short-run marginal asset --- primarily transformer --- degradation, and illustrates that voltage and/or ampacity congestion DLMC components are insufficient in providing the price signal that supports the system optimal network and DER schedules.

Third, it demonstrates that apart from optimal short-run scheduling, our approach can harvest otherwise idle DER reactive power compensation capabilities, increase distribution network DER (EV and PV) hosting capacity, and mitigate investments in distribution infrastructure that would be otherwise needed to support distribution utilities' obligation to serve. 

Finally, it provides numerical evidence on the significance of transformer degradation related intertemporally coupled DLMCs, suggesting that the introduction of intra-day markets may be worth considering.

The remainder of this paper is organized as follows.
Section \ref{DLMCs} discusses the use of DLMCs as price signals that capture the cost/benefit of intertemporal shifts in consuming/producing real/reactive power.
Section \ref{CaseStudy} presents the details of the case study considered in deriving numerical results and lists the alternative input data scenarios representing different EV/PV penetration levels.
In Section \ref{NumResults}, we present and discuss the main numerical-result-based findings, which we further enhance in Section \ref{FurtherNum}.
Section \ref{Conclusions} concludes and provides directions for future research.
In Appendix \ref{AppB}, we detail the sensitivity calculations employed in the DLMC unbundling.

\section{DLMCs as Price Signals} \label{DLMCs}

In this section, we employ DLMCs as price signals that provide DERs with sufficient information to self-schedule in a manner interpretable as a minimization (maximization) of their individual cost (benefit).
Indeed, DLMCs support optimal DER self-scheduling.
A point of caution worth making is that DLMCs are not necessarily prices that, when charged, will render the distribution network whole in terms of allowing it to recover its variable and fixed costs. 
The difference of DLMCs during periods $t$ and $t'$ simply represents the change in the system (marginal) cost if a unit of power consumption is transferred from period $t$ to period $t'$. 
For example, spatiotemporal DLMC-based prices may be set equal to a constant + the DLMC, as in a two-part tariff pricing approach.
As long as the price differences across time and location equal the DLMC and the value of the constant is selected to provide adequate total revenue, the two-part tariff design will also support optimal DER self-scheduling.

Let us consider rooftop solar PV $s$, which is offered a price $\hat \lambda_{t}^P$ for the provision of real power at time period $t$, denoted by $p_{s,t}$, and a price $\hat \lambda_{t}^Q$ for the provision of reactive power, denoted by $q_{s,t}$. 
The following optimization problem, referred to as \textbf{PV-opt}, describes a self-scheduling PV that adjusts its power factor to maximize its revenues from the provision of real and reactive power over the optimization horizon:\\ 
\textbf{PV-opt}:
\begin{equation} \label{PVopt}
\underset{p_{s,t},q_{s,t}}{\text{maximize}}  \sum_t { \left( \hat\lambda_{t}^{P} p_{s,t} + \hat\lambda_{t}^Q q_{s,t} \right)},
\end{equation}
\emph{subject to}:
PV constraints \cite[Eqs. (10)--(11)]{PartI}.

We note that PV-opt can be solved in parallel for each time period.
Evidently, if $\hat \lambda_{t}^P$ is negative, then the PV will select to produce no real power, setting $p_{s,t} = 0$.
If $\hat \lambda_{t}^Q$ is positive (negative), then the PV will adjust its power factor to provide (consume) reactive power, i.e., $q_{s,t} > 0$ ($q_{s,t} < 0$), so that the term $\hat\lambda_{t}^Q q_{s,t}$ is positive.

Similarly, let us consider EV $e$, which is offered a price $\hat \lambda_{t}^P$ for consuming real power at time period $t$, denoted by $p_{e,t}$, and a price $\hat \lambda_{t}^Q$ for consuming reactive power, $q_{s,t}$, at the node that it is connected during each time period.
The following optimization problem, referred to as \textbf{EV-opt}, describes a self-scheduling EV that employs its EV charger so as to minimize its net charging cost over the daily cycle:\\
\textbf{EV-opt}:
\begin{equation}
\underset{p_{e,t},q_{e,t}}{\text{minimize}} \sum_t  { \left( \hat \lambda_{t}^{P} p_{e,t} + \hat\lambda_{t}^Q q_{e,t} \right)},
\end{equation}
\emph{subject to}:
EV constraints \cite[Eqs. (12)--(17)]{PartI}.

Unlike PV-opt, EV-opt cannot be solved in parallel, as it includes intertemporal constraints.
Evidently, if $\hat \lambda_{t}^Q$ is positive (negative), then the EV will provide (consume) reactive power, i.e., $q_{e,t}$ will be negative (positive), so that the term $\hat\lambda_{t}^Q q_{e,t}$ is negative thus reducing the charging cost. 

Let us now assume that the solution of the operational planning problem described in Part I \cite{PartI} --- denoted by \textbf{Full-opt} in \cite[Subsection II-D]{PartI} --- is known, and the optimal DER schedules are $p_{s,t}^{*}, q_{s,t}^{*}$ for PVs and $p_{e,t}^{*}, q_{e,t}^{*}$ for EVs.
The DLMCs of the Full-opt solution at node $j$, time period $t$, reflecting the optimal DER schedules, are $\lambda_{j,t}^{P*}$ and $\lambda_{j,t}^{Q*}$.
It is easy to show that if $\lambda_{j,t}^{P*}$ and $\lambda_{j,t}^{Q*}$ were the prices announced to each PV $s$ and EV $e$, then the optimal PV and EV schedules of Full-opt would also be optimal for the PV-opt and EV-opt problems, respectively.
To verify the above, consider the partial Lagrangian of Full-opt developed in \cite{PartI}.
It is obtained by appending the real and reactive power balance constraints \cite[Eqs. (2)--(3)]{PartI} and substituting the net demand variables with \cite[Eqs. (8)--(9)]{PartI}.
The terms that include the DER (PV/EV) variables in the partial Lagrangian are as follows:
\begin{equation*} 
\sum_{e,t} { \left( \lambda_{j_e,t}^{P} p_{e,t} + \lambda_{j_e,t}^Q q_{e,t} \right)} 
- \sum_{s,t} { \left( \lambda_{j_s,t}^{P} p_{s,t} + \lambda_{j_s,t}^Q q_{s,t} \right)},
\end{equation*}
where $j_s$ refers to the node that PV $s$ is installed, and $j_e$ to the node that EV $e$ is connected at time period $t$.
The first sum coincides with the objective function of EV-opt, for $\lambda_{j_e,t}^{P} =  \hat \lambda_{t}^{P}$ and $\lambda_{j_e,t}^{Q} = \hat \lambda_{t}^{Q}$, summed for all EVs.
Similarly, the second sum coincides with the objective function of PV-opt summed for all PVs;
the minus sign is obtained if we convert PV-opt to a minimization problem. 
Also, we note that PV constraints appear in both PV-opt and Full-opt, and that EV constraints appear in both EV-opt and Full-opt.

Let us consider the optimality conditions of Full-opt that involve DLMCs and DER variables (it is convenient to consider the partial Lagrangian representation), and the respective optimality conditions of PV-opt and EV-opt.
It is rather straightforward to see that if we replace in PV-opt and EV-opt parameters $\hat \lambda_{t}^{P}$ and $\hat \lambda_{t}^{Q}$ with the DLMCs derived by the optimal solution of Full-opt, say $\lambda_{j,t}^{P*}$ and $\lambda_{j,t}^{Q*}$, with $j$ referring to appropriate nodes $j_s$ and $j_e$, respectively, then the optimality conditions of PV-opt and EV-opt will also be encountered in the optimality conditions of Full-opt that refer to PV and EV variables, respectively.
Hence, an optimal solution of the Full-opt problem will also be optimal for the PV-opt and EV-opt problems if the prices announced to PVs and EVs reflect the respective DLMCs of the Full-opt solution.
   
Although only PVs and EVs are modeled in this paper as representative DER examples, other DERs can be treated similarly.
The individual optimization problems suggest the following interpretation: 
Nodal marginal costs can be construed as prices that elicit a price-taking DER to adapt fully and self-schedule to its socially optimal real/reactive power profile. 
In other words, if we were able to determine these socially optimal spatiotemporal DLMCs and charge DERs on DLMC-based-prices, DERs would self-schedule in a manner that is optimal for the system as a whole.

We wish to note that DLMC-based pricing has attracted strong criticism that points to undesirable DLMC volatility. 
However, our numerical results provide overwhelming evidence that DER schedule adaptation to DLMC-based prices removes the volatility, which is indeed observed only when DERs schedule themselves in a manner that is not adaptive to the spatiotemporal DLMCs.

\section{Case Study} \label{CaseStudy}
\begin{figure*}[ht]
\centering
\includegraphics[width=6.5in]{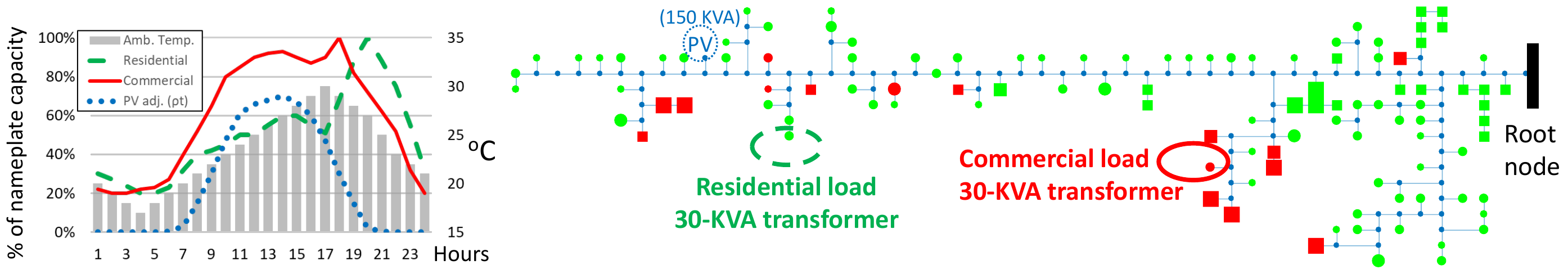}
\caption{Feeder topology diagram (307 nodes). 110 transformers: round (pole); square (pad); green (residential); red (commercial); sizes reflect nameplate capacity. Commercial and residential load profiles, PV adjustment factor ($\rho_t$), and ambient temperature $ \theta^H_t$. PV installation 150 KVA (p.f. = 1).} 
\label{figFeeder}
\end{figure*}
\begin{table}[tb] 
	\caption{Aggregate Line and Transformer Data}  \label{tab1} 
\centering
\includegraphics[width=3.45in]{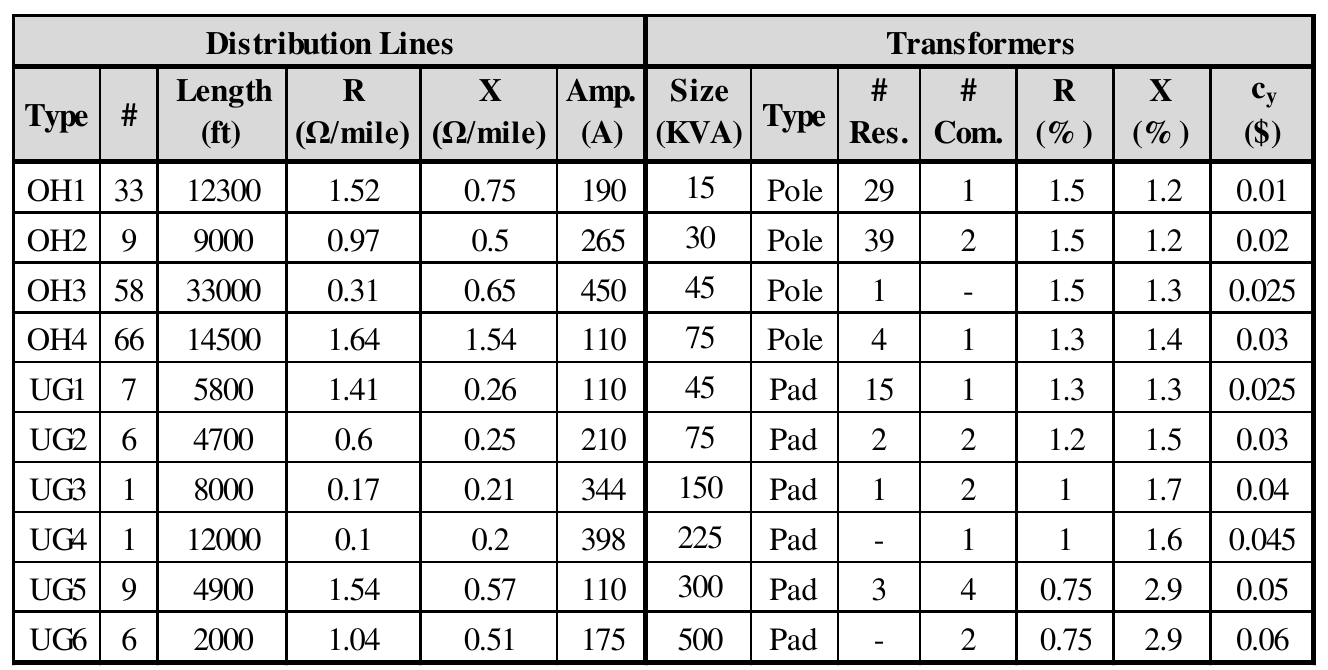}
\end{table}

In this section, we illustrate the application of the proposed model on an actual 13.8 KV feeder of Holyoke Gas and Electric (HGE), a municipal distribution utility in MA, US. 
In Subsection \ref{FeederData}, we list the input data of the feeder. 
In Subsection \ref{Options}, we introduce several reasonable, though sub-optimal, EV and PV scheduling options that we compare to the optimal schedule. 

\subsection{Input Data} \label{FeederData}

The feeder topology is shown in Fig. \ref{figFeeder}.
Aggregate line and transformer data are listed in Table \ref{tab1}. 
Lower and upper voltage limits are set at 0.95 and 1.05 p.u., respectively.
Fig. \ref{figFeeder} also shows the commercial and residential load profiles (as a percentage of the transformer nameplate capacity), the PV adjustment factor, $\rho_t$, and the ambient temperature. 
The load profiles are obtained from \cite{CIGRE}, assuming a 0.95 (0.85) power factor for the residential (commercial) node.
LMPs range from 25.59 to 53.48 (\$/MWh).
The opportunity cost for reactive power is assumed at 10\% the value of the LMP.

Numerical experimentation focuses on two selected nodes depicted in Fig. \ref{figFeeder}, representing two 30-KVA transformers serving commercial and residential loads.
For both transformers, the ratio of load losses at rated load to no load losses is $R=5$, the rise of the top-oil temperature over ambient temperature at rated load is $\Delta \bar \theta^{TO} = 55$, and the rise of HST over top-oil at rated load is $\Delta \bar \theta^{H} = 25$.

To make results easier to follow while emphasizing the local effect of EVs and PVs on a distribution feeder, we consider different levels of EVs and PVs connected exclusively to these two nodes.
At the commercial node, EVs are connected 9am--5pm and require charging 12 KWh.
At the residential node, EVs are connected 7pm--7am and require charging 18 KWh. At the time of departure, EVs must be fully charged.
EV battery capacity is 24 KWh, the maximum charging rate 3.3 KW/h, and the charger capacity 6.6 KVA.
PVs are assumed to be 10 KVA rooftop solar.

To model initial conditions in the daily cycle reasonably, transformer temperatures were required to coincide at the beginning and at the end of the cycle, $t=0$, and $t=24$; a similar constraint was imposed on EV battery State of Charge.

\subsection{Scheduling Options} \label{Options}

We consider 4 EV/PV scheduling alternatives:

(1) ``\textbf{BaU}'' (Business as Usual): EVs ``dumb'' charge at full rate upon arrival with unity power factor (p.f. = 1). 
PVs operate with p.f. = 1.

(2) ``\textbf{ToU}'' (Time-of-Use): EVs charge to minimize real power cost relative to the selected LMP trajectory while maintaining a p.f. = 1. PVs operate with p.f. = 1.

(3) ``\textbf{PQ-opt}'': EVs/PVs are scheduled in order to minimize real and reactive power cost, subject to voltage and ampacity constraints while ignoring transformer degradation cost.
PQ-opt can be thought of as a traditional line loss minimization problem akin to the day-ahead problem presented in \cite{BaiEtAl-DLMP}.
It can be alternatively thought of as a special case of Full-opt below without transformer degradation costs and related constraints.

(4) ``\textbf{Full-opt}'': EVs/PVs are scheduled by adapting to full DLMCs, namely including transformer degradation costs.
This adaptation is given by the solution of Full-opt detailed in \cite{PartI}.

For each of the above scheduling options, we analyze and compare several DER penetration scenarios at the nodes of interest:
EVs equal 0, 3, and 6, and rooftop PV installations equal 0, 30, and 60 KVA (i.e., 0, 3, and 6 units of 10-KVA rooftop solar, respectively).
The scenario with 0 EVs and PVs serves as the base case for reference purposes.

Section \ref{NumResults} reports and compares numerical results across scenarios and scheduling options.
Additional numerical results focusing on the intertemporal impact of the transformer DLMC component, ampacity congestion, and higher DER (EV/PV) penetrations of up to 12 EVs and 120 KVA PV are presented in Section \ref{FurtherNum}.

\section{Numerical Results} \label{NumResults}
\begin{table*}[ht] 
	\caption{Aggregate System Cost Difference (in \$) and LoL (in hours)}  \label{tab3a} 
\centering
\includegraphics[width=6.5in]{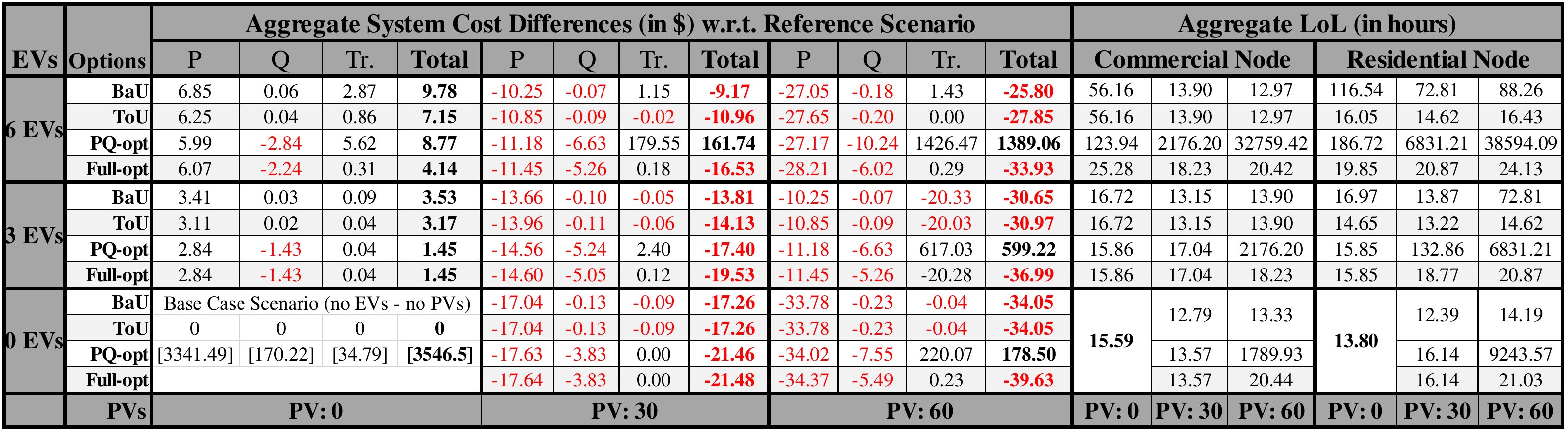}
\end{table*} 
DER schedules were obtained by solving the corresponding optimization problem on a Dell Intel Core i7-5500U @2.4 GHz with 8 GB RAM, using CPLEX 12.7. Solution times were up to 10 sec.
Although the convex relaxation of the branch flow model was exact in all instances, we refer the interested reader to recent works that propose remedies when the relaxation is not exact \cite{Abdelouadoud_EtAl_2015, HuangEtAl_2017, WeiEtAl_2017, NickEtAl_2018, YuanEtAl-DLMP}.
Sensitivities that unbundle DLMCs to additive components were obtained by solving the linear systems detailed in Appendix \ref{AppB}.

Table \ref{tab3a} shows the system-wide daily cost differences (in \$) for real/reactive power (P/Q), transformer degradation, and total cost relative to the base case (no EVs and no PVs).
For readability, negative values are shown in red.
Table \ref{tab3a} reports also the aggregate daily Loss of Life (LoL) of the two 30-KVA transformers (commercial and residential) under different scheduling options and EV/PV penetrations.

As expected, Full-opt achieves the lowest total cost, since it co-optimizes real/reactive power and transformer cost and maintains low LoL.
PQ-opt achieves the lowest combined real and reactive power costs (P and Q columns in Table \ref{tab3a}), but exhibits some very high values of aggregate LoL and transformer degradation cost, as a result of high reactive power provision.
Unlike PQ-opt and Full-opt, ToU and BaU do not take advantage of DER reactive power provision capabilities, but in some cases, particularly in the residential node, achieve low aggregate LoL.
 
In the commercial node, BaU and ToU produce identical EV schedules (since by coincidence LMPs are increasing during the day).
In the residential node, ToU exhibits a 3-hour shift relative to the BaU EV profile (EVs start charging at 10pm). 
For the specific LMPs and EV/PV penetrations considered, the ToU schedule results in sustainable transformer LoL performance;
different LMPs though and/or higher penetration levels may result in worse performance of the ToU option. 
Since we focus on understanding the granularity of DLMCs, we did not consider multiple EV profiles and arrival times that would demonstrate the potential synchronization under the ToU schedule (see e.g., the simulations in \cite{HilsheyEtAl_2013}) resulting in worse performance relative to the BaU schedule.
For the high EV/PV penetration presented in Section \ref{FurtherNum}, however, high LoL is associated with ToU.

We focus next on real and reactive power DLMCs.
It is important to note that with the exception of Full-opt generated schedules, all other DER schedules are not adapted to DLMCs.
For BaU, ToU, and PQ-opt DER schedules, DLMCs simply represent the marginal cost of delivering an incremental unit of real and reactive power to/at a specific network location and time, conditional upon the given DER (EV/PV) schedule. 
As such, DLMCs under all the non-adaptive schedules should be interpreted as ex post spatiotemporal marginal costs indicating desirable directions of change.
So for a given DER schedule (BU, ToU, PQ-opt) obtained according to Subsection \ref{Options}, DLMCs are calculated solving a partial Full-opt problem with given/fixed EV/PV schedules.
As such the DLMCs can be interpreted as economically/socially efficient financial incentives indicating a direction for improving system costs.
For example, a DLMC at a specific node and time which is higher than the Full-opt DLMC signals indicates that less power should be consumed then and there, if at all possible.

We next proceed to investigate and discuss EV-only (no PVs) scenarios in Subsection \ref{EVonly}, PV-only (no EVs) scenarios in Subsection \ref{PVonly}, and EV-PV synergy scenarios in Subsection \ref{EVPVSynergy}. 

\subsection{EV Only} \label{EVonly}
\begin{figure}[tb]
\centering
\includegraphics[width=3.2in]{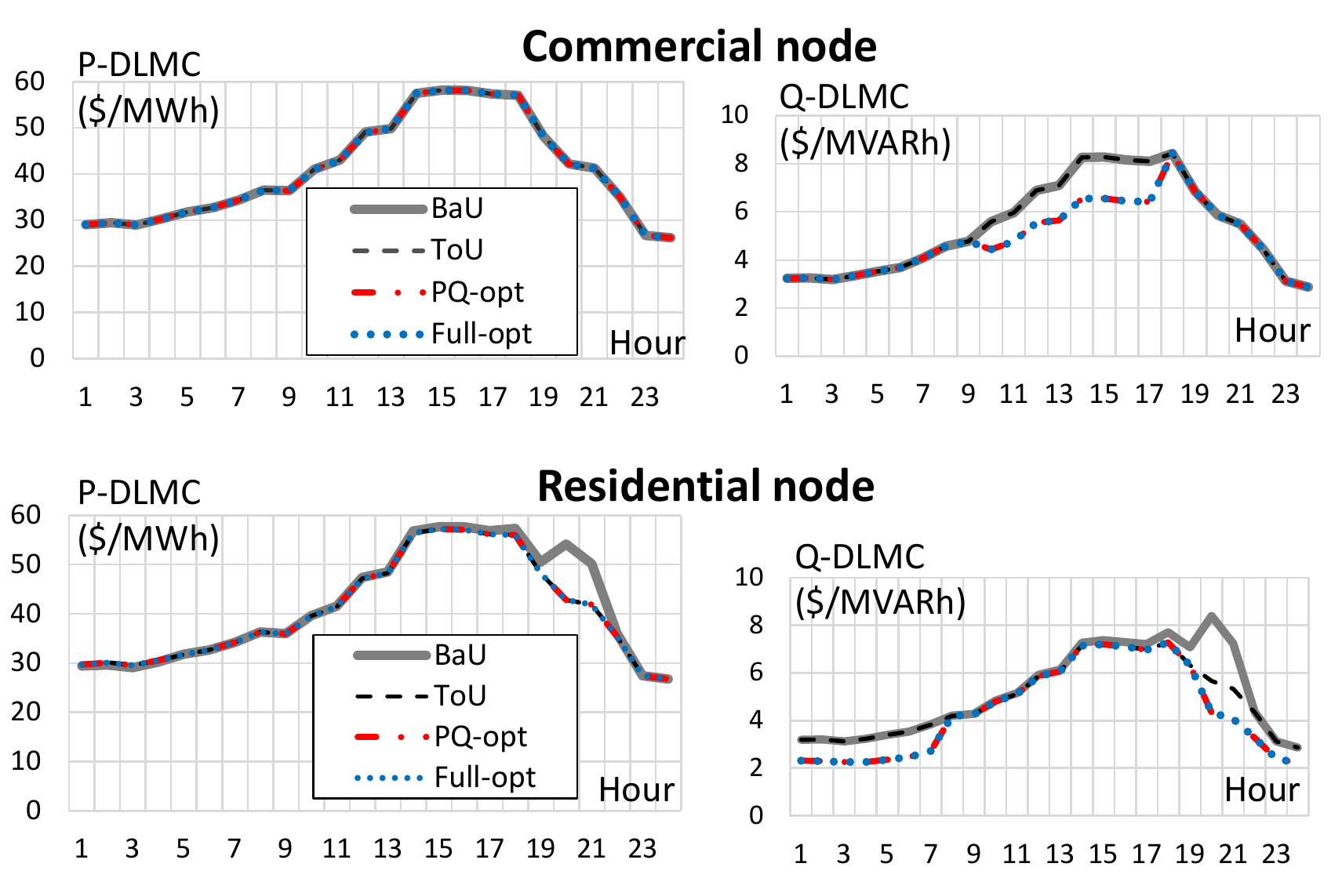}
\caption{DLMCs, commercial and residential nodes, 3 EVs.}
\label{figEV3A}
\end{figure}

In Fig. \ref{figEV3A}, we show P-DLMCs and Q-DLMCs for a relatively low 3-EV penetration scenario. 
The differences are rather small and are related to difference in the provision of reactive power.
At the commercial node, P-DLMCs are practically the same under all scheduling options. 
Q-DLMCs are higher under BaU/ToU compared to PQ-opt/Full-opt, indicating that higher reactive power provision would be desirable.
At the residential node, P-DLMCs associated with BaU EV schedules are higher for hours 20--22, pointing to the desirability of a 3-hour shift in EV consumption;
this coincides with the result of other scheduling options (EVs start charging at 10pm).
Q-DLMCs associated with BaU/ToU EV schedules suggest the desirability of providing reactive power;
this coincides with the results of Full-opt and PQ-opt that indeed provide reactive power and decrease system-wide cost.

\begin{figure}[tb]
\centering
\includegraphics[width=3.2in]{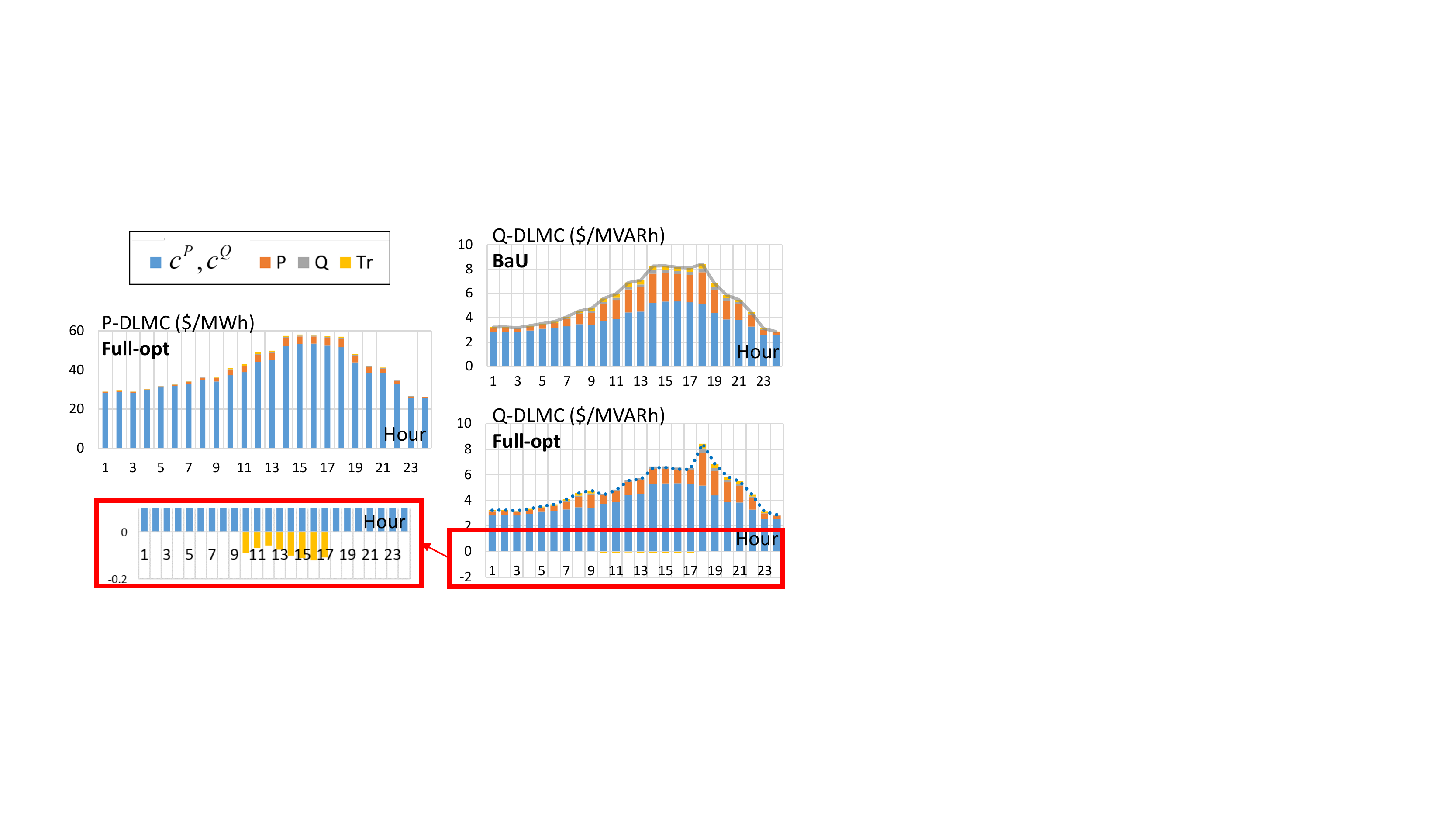}
\caption{DLMC components, commercial node, 3 EVs.}
\label{figEV3B}
\end{figure}

Fig. \ref{figEV3B}, reports the marginal cost components of P-DLMCs for Full-opt (since they are similar across all scheduling options), and Q-DLMCs for BaU (same as ToU) and Full-opt (same as PQ-opt), at the commercial node.
P-DLMC components range from 1.7\% to 7.5\% for real power losses, from 0.2\% to 0.9\% for reactive power losses and from 0.2\% to 1.7\% for transformer degradation. 
For Q-DLMCs, we observe that during hours 10--17 (when EVs are plugged in), the marginal cost components for real and reactive power are lower under Full-opt, while the transformer degradation component becomes negative.
This is interestingly associated with a reverse reactive power flow at that node and those time periods rendering the sensitivity of the current \emph{w.r.t.} reactive power net demand negative. 
Q-DLMC components under Full-opt range from 8\% to 30.7\% for real power losses, from 1\% to 3.7\% for reactive power losses, and from -2\% to 5\% for transformer degradation.

\begin{figure}[tb]
\centering
\includegraphics[width=3.2in]{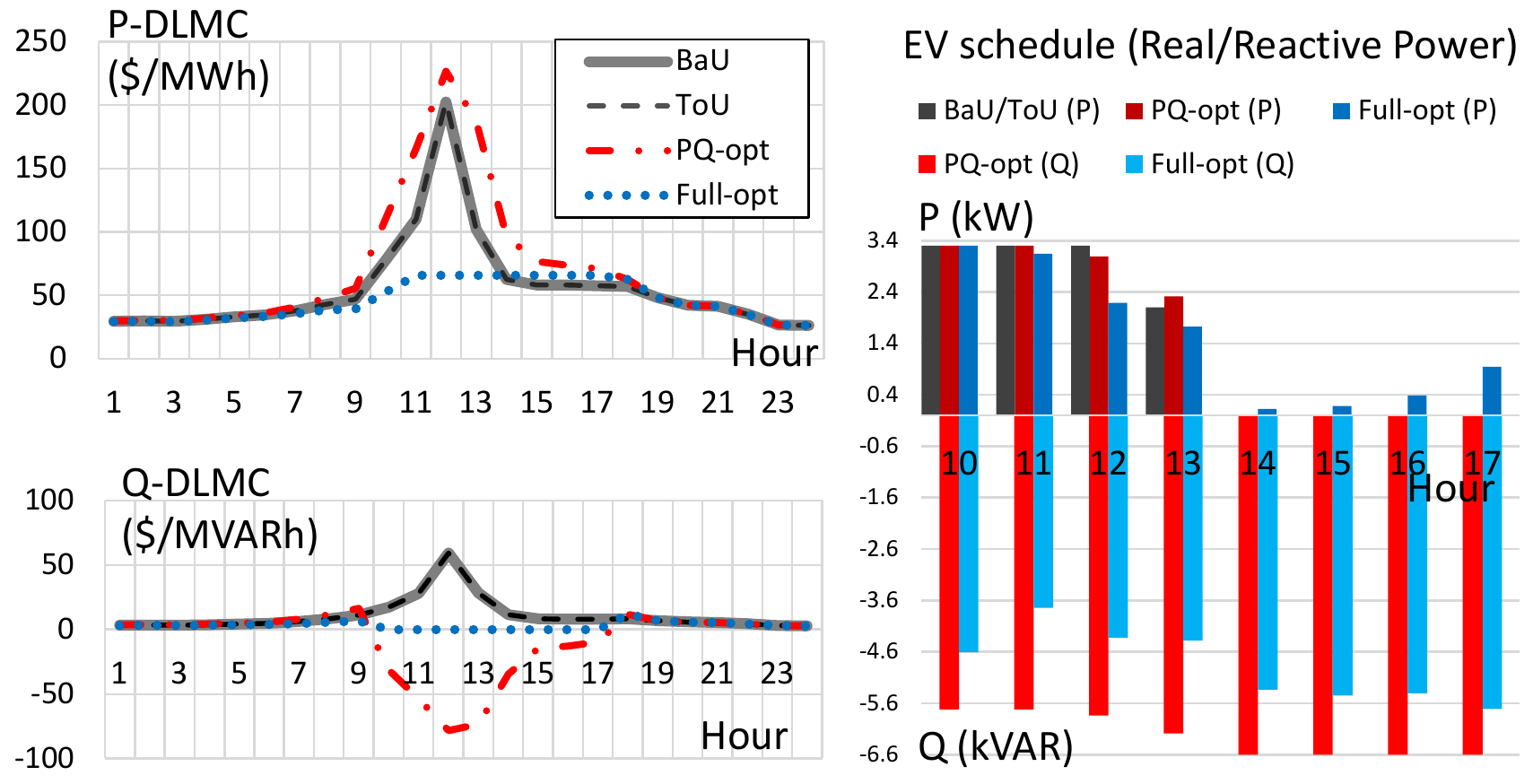}
\caption{DLMCs at commercial node, 6 EVs.}
\label{figEV6C}
\end{figure}

EV penetration is more impactful under the 6-EV penetration scenario. PQ-opt is associated with a 5 to 9 fold higher aggregate transformer LoL relative to Full-opt (see Table \ref{tab3a}).
Fig. \ref{figEV6C} reports the 6-EV penetration DLMCs at the commercial node, and their schedule during hours 10--17 that they are plugged in.
P-DLMCs exhibit spikes under BaU/ToU and PQ-opt implying that charging rates are too high, whereas Full-opt exhibits a smoother profile.
Q-DLMCs exhibit positive spikes under BaU/ToU, implying that EVs should provide reactive power. 
Negative spikes under PQ-opt suggest that EVs provide excessive reactive power. 
We note that Full-opt results in practically zero Q-DLMCs ($\to 0^+$) while EVs are plugged in, thus implying high reactive power availability. 
Nevertheless, despite the fact that the imputed reactive power income --- 
recall that EV-opt represents self-scheduling by an EV that minimizes its net charging cost --- 
decreases with tanking Q-DLMCs, Full-opt results in lower P-DLMCs as well; the net impact on the overall system cost from EV charging is lower and EVs realize less imputed income from Reactive Power ``sales'' but overall their effective charging costs are lower!
Full-opt Q-DLMCs support the system-optimal solution by incentivizing EVs to provide reactive power at a rate that is lower than their charger capability,
whereas PQ-opt ``brute-forces'' EVs to fully utilize their inverter capability resulting in negative Q-DLMCs that significantly increase the imputed cost of EVs 
(negative Q-DLMCs combined with negative reactive power increase the cost for EVs in EV-opt).

\begin{figure}[tb]
\centering
\includegraphics[width=3.2in]{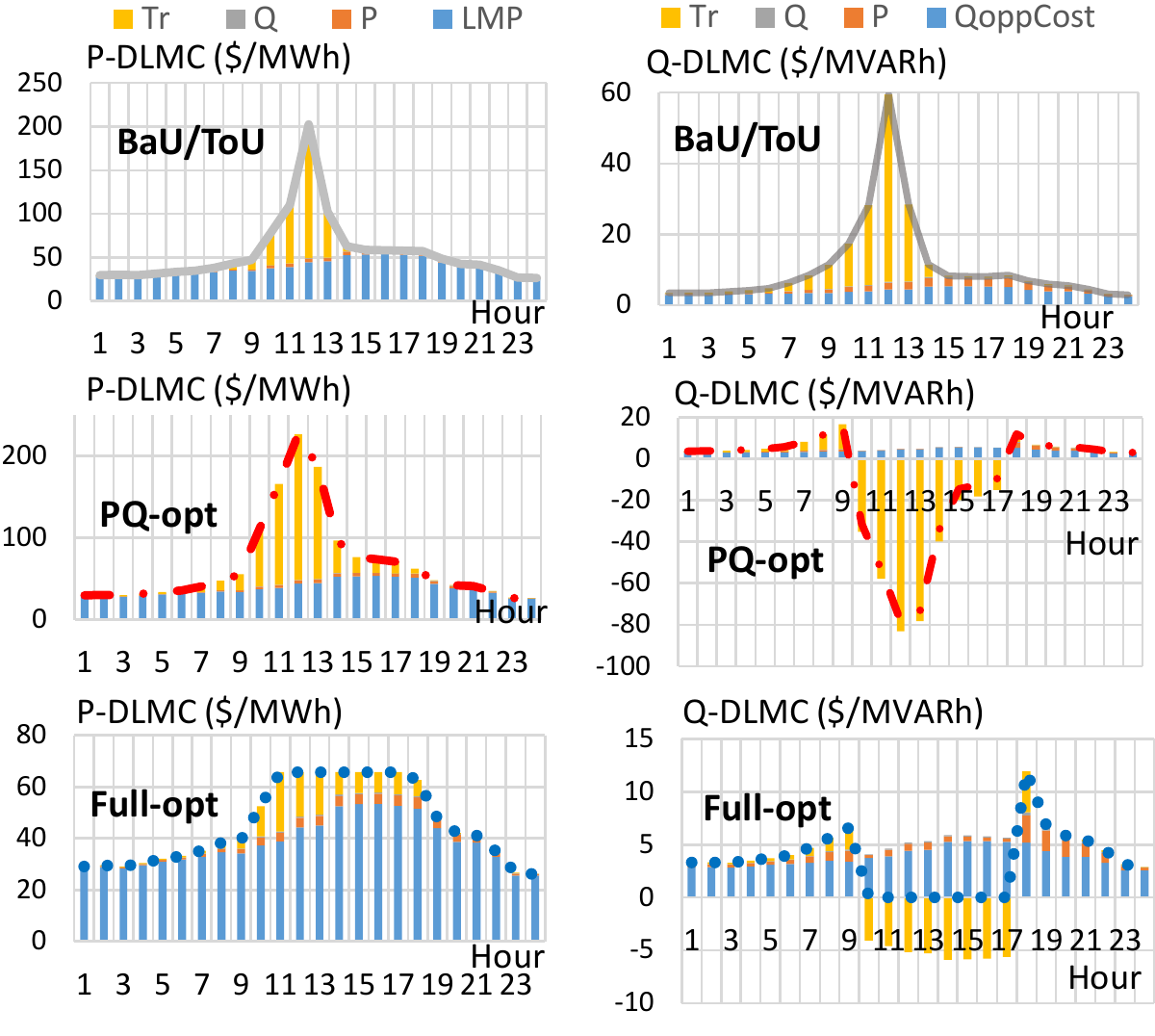}
\caption{DLMC components, commercial node, 6 EVs.}
\label{figEV6D}
\end{figure}

Fig. \ref{figEV6D} shows the DLMC components at the commercial node again under a 6-EV penetration.
It is noteworthy that the marginal transformer degradation cost varies across scenarios and scheduling options illustrating the incentives implied by the respective DLMCs.
Indeed, the spikes (positive/negative) are caused by the transformer degradation component.
In general, low transformer costs incurred by Full-opt across all hours are associated with smoother DLMC profiles.

Notably, although voltage and ampacity DLMC components are zero, i.e., the respective constraints are not binding, transformer degradation is high. 
This is an interesting observation showing that traditional line loss minimization approaches employing the Volt/VAR control capabilities of smart inverters and capacitors, such as in PQ-opt, may ``mask'' the significant impact of transformer degradation on DLMCs, and impose potentially excessive cost on the distribution network assets.

Similar remarks can be made for the residential node. 
They are omitted due to space considerations.

\subsection{PV Only} \label{PVonly}

The 30-KVA PV penetration scenario has rather insignificant impact on network costs (see Table \ref{tab3a}).
The 60-KVA penetration scenario, on the other hand, is quite impactful, with PQ-opt yielding very high aggregate LoL values (see Table \ref{tab3a}) compared to BaU or Full-opt. 
Such high LoL values imply excessive transformer overloading;
in this case, imposing ampacity constraints --- even below overcurrent protection --- would worth considering.
We elaborate further in Section \ref{FurtherNum}.

\begin{figure}[tb]
\centering
\includegraphics[width=3.4in]{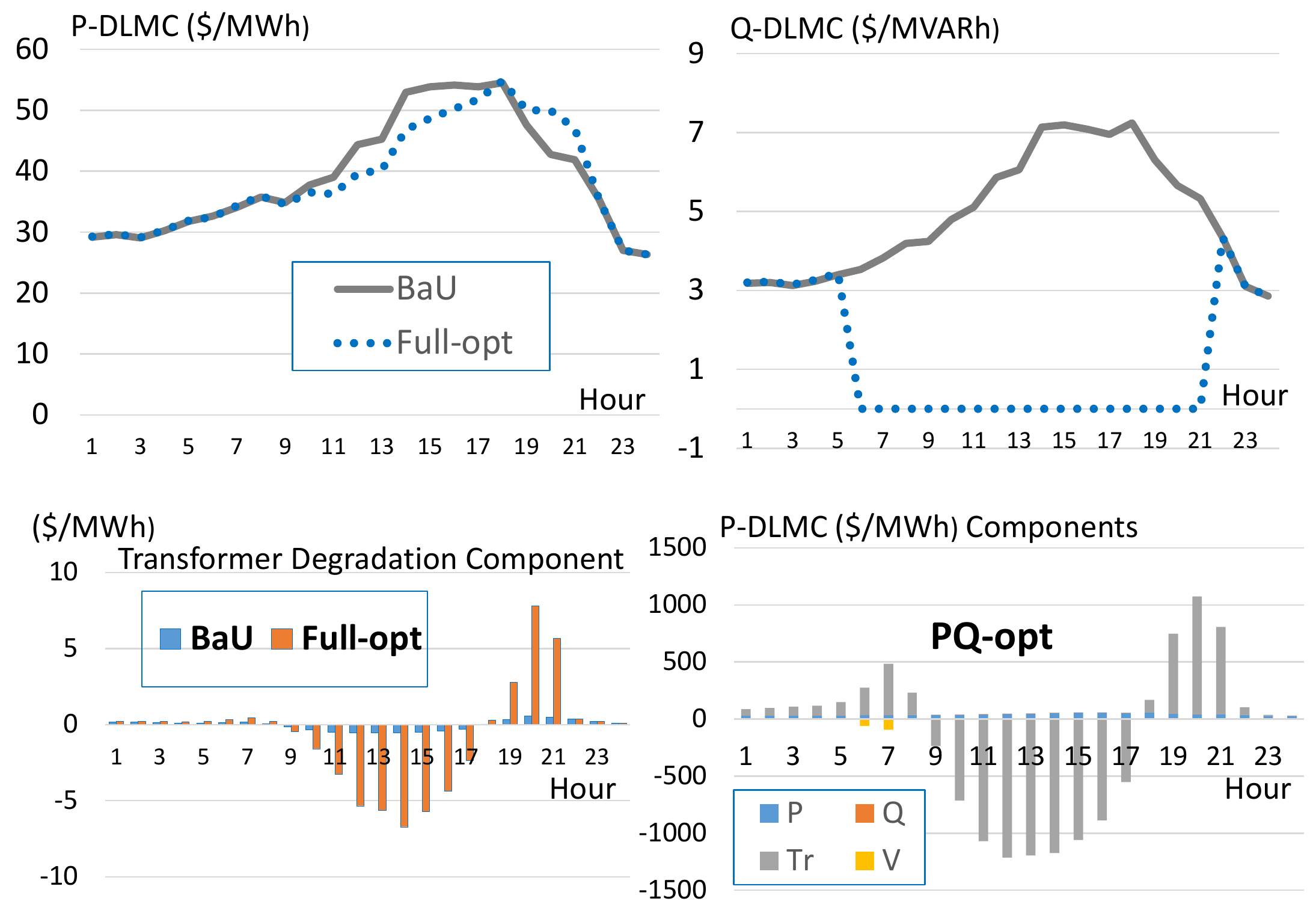}
\caption{DLMCs, residential node, 60 KVA PV.}
\label{figPV60}
\end{figure}

In Fig. \ref{figPV60}, we illustrate DLMCs at the residential node for the 60-KVA PV penetration.
We present P-DLMCs and Q-DLMCs under BaU and Full-opt (top figure), the P-DLMC transformer degradation component for BaU and Full-opt (bottom left), and the P-DLMC components for PQ-opt (bottom right).

We observe that BaU and Full-opt result in similar P-DLMCs (top left); the differences are explained by the transformer degradation component (bottom left).
For BaU and Full-opt, the transformer component becomes negative during hours 9--17, when there is a reverse real power flow (PV generation exceeds the load) and the sensitivity of the current \emph{w.r.t.} net demand becomes negative.
The transformer component under PQ-opt (bottom right) is ``amplified'' --- at least two orders of magnitude higher than Full-opt.
We also observe a negative voltage component under PQ-opt at hours 6 and 7, when the upper bound of the voltage is binding, indicating a negative sensitivity of the voltage \emph{w.r.t.} net demand.
Indeed, increasing net demand would result in a voltage drop.
However, the transformer degradation component at those hours is positive, indicating a positive sensitivity of the current \emph{w.r.t.} net demand, and is about 4 to 5 times larger in magnitude than the voltage component.

Similar to the 6-EV scenario, DLMC differences are caused by reactive power provision.
Positive Q-DLMCs under BaU (top right) encourage the provision of reactive power.
Under PQ-opt, Q-DLMCs (not shown in Fig. \ref{figPV60}) become highly negative, suggesting that reactive power provision is excessive (PVs fully utilize their inverter capabilities).
Under Full-opt, PVs provide reactive power as needed to minimize system cost; 
Q-DLMCs tank, indicating that there is excess installed capacity and it is optimal to underutilize it.
Full utilization, as is the case under the PQ-opt schedule, renders negative Q-DLMCs and interacts with P-DLMCs to result in a large transformer degradation component.

\subsection{EV-PV Synergy} \label{EVPVSynergy}

With the exception of PQ-opt at the residential node, the 3-EV and 30-KVA PV scenarios are easily sustainable with low aggregate LoL. 
The EV-PV synergy at hours 6, 7, 20, and 21 results in high reverse reactive power flows that overload the transformer. 
P-DLMCs exhibit positive spikes and Q-DLMCs negative spikes at those hours. 
Increasing EVs to 6 and/or PVs to 60 KVA results in very high LoL under PQ-opt scheduling.
On the other hand, as expected, Full-opt achieves sustainability with low LoL under all scenarios.

\section{Further Numerical Investigation} \label{FurtherNum}

In this section, we investigate further the intertemporal behavior of the transformer degradation marginal cost component (Subsection \ref{intertemp}), transformer ampacity congestion (Subsection \ref{ampcong}), and
increased EV/PV penetration (Subsection \ref{increased}).

\subsection{Transformer Degradation Intertemporal Impact} \label{intertemp}

In general, the transformer degradation component at a specific location includes the impact on all transformers through weights that are proportional to the sensitivities of the current \emph{w.r.t.} net demand.
For ease of exposition, in our numerical experimentation, we consider the impact only on the transformer located at a specific node.\footnote{
This is not an unreasonable approximation, since in our test case, overloads are imposed on two rather distant transformers.}
Due to space considerations, we report on the P-DLMC component;
similar analysis holds for the Q-DLMC component.

Let us consider transformer $y$ denoted by line $ij'$. 
The transformer degradation component (including only the impact on transformer $y$) at node $j'$, time period (hour) $t'$ is given by:
\begin{equation} \label{Margfyt1}
\begin{split}
 \pi_{y,t'} & \frac{\partial l_{y,t'}}{\partial p_{j',t'}}
 = \rho_y \left[ \frac{55}{6} { \left( \frac{3}{4} \right)}^{T-t'}  
 \right] \frac{\frac{\partial l_{y,t'}}{\partial p_{j',t'}} }{l_y^N} \\
 +    c_y& 
 \left[
 \frac{55}{6}  \tilde \alpha_{y,t'} 
 + 20 \tilde \alpha_{y,t'}
 + \sum_{\tau = 1}^{T-t'} { \left( \frac{3}{4} \right)}^{\tau} \tilde \alpha_{y,t'+ \tau}
 \right] 
 \frac{\frac{\partial l_{y,t'}}{\partial p_{j',t'}} }{l_y^N},
 \end{split}
\end{equation}
where $\rho_y$ is the dual variable of the constraint requiring the HST at the end of hour 24 to equal the HST at the beginning of hour 1, and as such captures the impact of increasing the HST in hour 24 on future -- next day -- HST and transformer degradation. 
The term $\tilde \alpha_{y,t} = \frac{1}{c_y}\sum_{\kappa} \xi_{y,t,\kappa} \alpha_{\kappa}$ can be viewed as an adjusted slope of the aging acceleration factor,\footnote{
Recall that $\sum_{\kappa} \xi_{y,t,\kappa} =c_y$. 
In the example used in \cite{PartI}, the slopes $\alpha_{\kappa}$ range from 0.009 for $\kappa = 1$ to 22.37 for $\kappa = 8$.}
and $l_y^N$ is the nominal current (squared).

First, we note that the component depends on the sensitivity of the current \emph{w.r.t.} net demand at hour $t'$, which determines the sign of the component, normalized by the transformer nominal current.
Second, we note that the component depends on the slope of the aging acceleration factor at hour $t'$, as well as the respective slopes of subsequent hours.
The impact that is related to hour $t'$ includes the contribution of the top-oil temperature and the winding.
In fact, there is a fixed ratio between the two that depends on the problem/transformer parameters, which in our case is equal to 2.18 (the contribution of the winding is higher).
The impact that is related to subsequent periods refers only to the top-oil temperature dynamics.
Interestingly, although this impact decays with ${ \left( \frac{3}{4} \right)}^{\tau}$, for each hour $\tau$ that follows hour $t'$, it can become significant whenever the slopes $\tilde \alpha_{y,t'+ \tau}$ become quite large.
As we show next, the impact that extends beyond the optimization horizon that is related to dual variable $\rho_y$ may also become significant.

\begin{figure}[tb]
\centering
\includegraphics[width=3.4in]{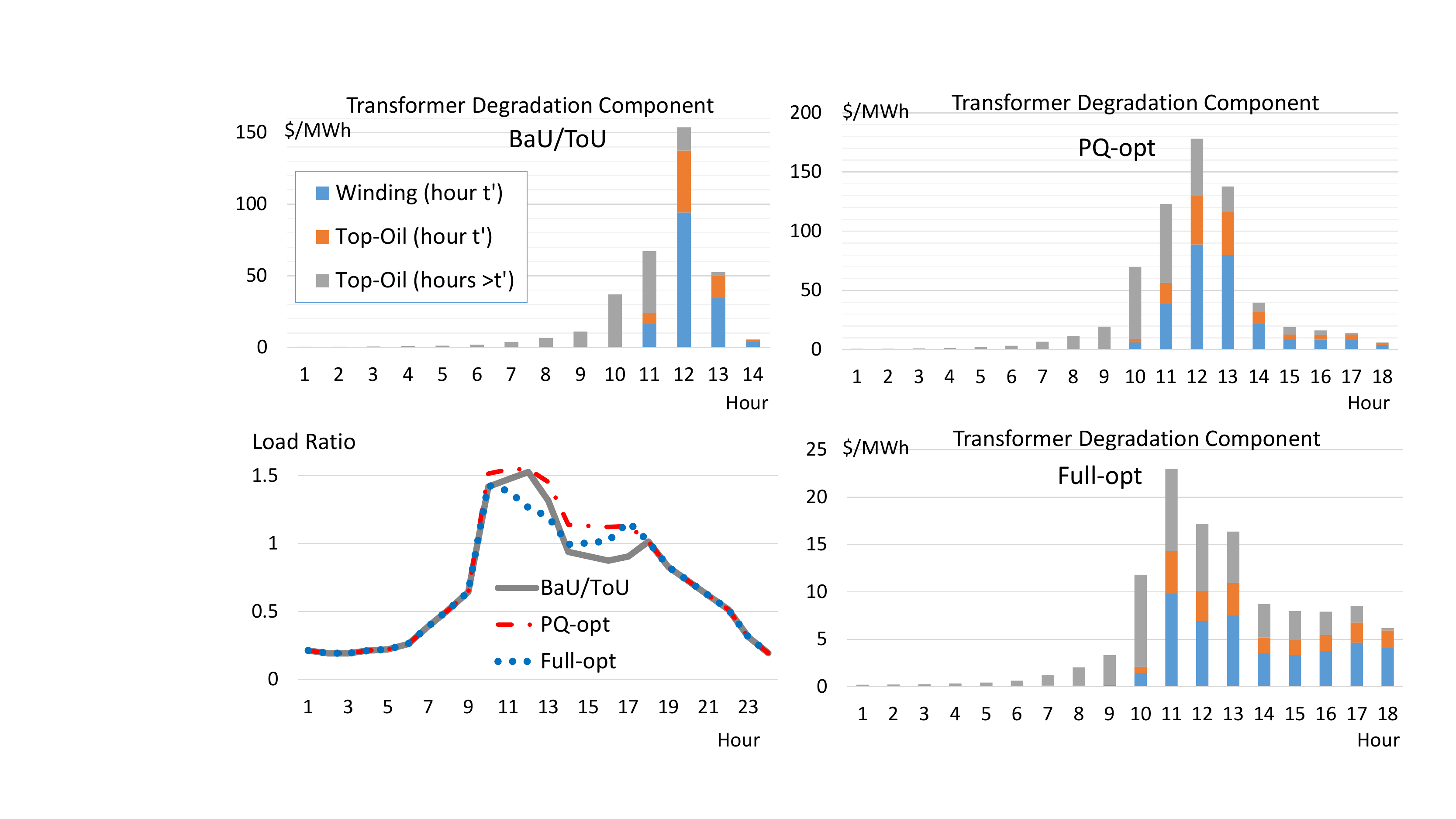}
\caption{P-DLMC transformer degradation components, commercial node, 6 EVs. Load Ratio $l_{y,t}/l_y^N$.}
\label{figEV6E}
\end{figure}
Let us consider the P-DLMCs for the 6-EV penetration scenario, at the commercial node, shown in Fig. \ref{figEV6D} (left).
For this scenario, we show in Fig. \ref{figEV6E} the load ratio and the sub-components of the transformer degradation.
An interesting remark is that the transformer component becomes significant before the EVs plug in, e.g., at hour 9, due to the impact of the top-oil temperature on subsequent hours. 

Under BaU/ToU schedules, the transformer component becomes 10.96 (\$/MWh) at hour 9, accounting for 23.3\% of the P-DLMC.
The component is still significant even at hours 7 and 8, accounting for 9.95\% and 15.3\% of the P-DLMC, respectively.
It becomes as high as 153.75 (\$/MWh) at hour 9, and diminishes after hour 14 as the transformer cools down.
At hours that precede the overloading hours (due to EVs) and at the first two hours of charging, the impact on the subsequent hours dominates the transformer component.
The winding and top-oil impact become dominant at hour 12.

PQ-opt exhibits similar results which are amplified and extended (in time duration) by the excessive brute-force-driven reactive power provision.
The transformer component becomes significant starting at hour 6 (9.3\% of the P-DLMC), it reaches a maximum of 178.1 (\$/MWh) at hour 12, the impact of subsequent hours remains large until hour 11, and then diminishes after hour 18.

Similar intertemporal effects are observed under the Full-opt, although the behavior is much smoother compared to PQ-opt.
The impact on subsequent hours makes up the major part of the P-DLMC until hour 10, and even though the transformer component becomes higher relative to BaU/ToU at hours 14--18 (ranging from 3 to 8 \$/MWh), the profile is smoother.

\begin{figure}[tb]
\centering
\includegraphics[width=3.2in]{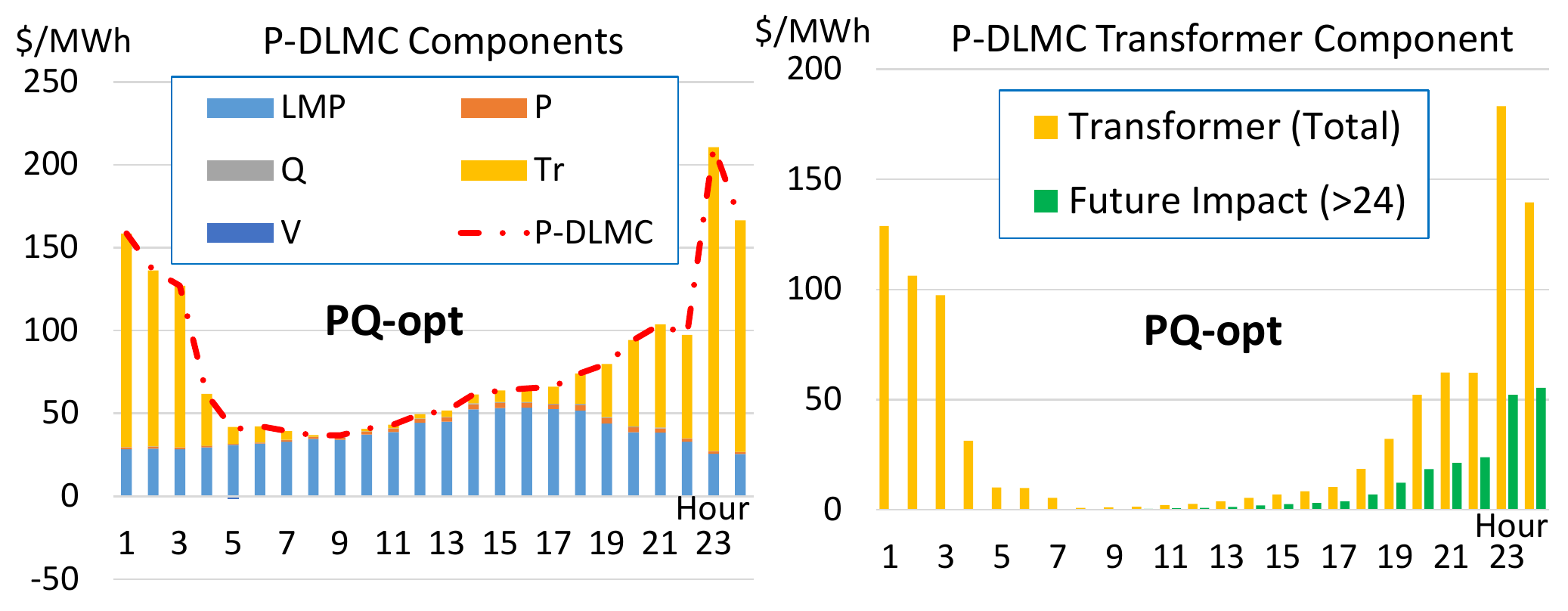}
\caption{P-DLMC components and transformer degradation future impact, residential node, 6 EVs.}
\label{figFuture}
\end{figure}

Fig. \ref{figFuture} shows the residential node P-DLMC components (left) for 6-EVs scheduled under PQ-opt. 
Interestingly, the impact on the early-next-day future transformer degradation component (right), of the suboptimal EV schedule during the late hours of the day is significant --- it becomes increasingly significant as we approach the end of the horizon.
It explains up to 33\% of the P-DLMC at hour 24, and less at previous hours (e.g., 20\% at hour 21, 15\% at hour 19, 5\% at hour 15).


The intertemporal impact of LoL and particularly the impact that extends beyond the end of a daily cycle on the DLMCs during the 24-hour daily cycle argues in favor of considering the introduction of intra-day markets.

\subsection{Transformer Ampacity Congestion} \label{ampcong}

In Section \ref{NumResults}, we did not impose hard transformer ampacity constraints, which may have been the culprit for some very high LoL values under PQ-opt schedules.
Hence, we are revisiting the issue as, in practice, excessive transformer loading may be disallowed by overcurrent protection, and we repeat our DLMC analysis by applying a transformer ampacity constraint equal to twice the value of the nominal current.
In most cases, the aggregate PQ-opt results are improved, but not enough to approach Full-opt behavior.
Indicatively, the 6-EV and 30-KVA PV penetration scenario at the residential node yields aggregate LoL of 4589 hours (compared to 6831 hours of the unconstrained case). 
The ampacity constraint is binding for 4 hours, with the ampacity congestion component in the same ``direction'' (sign) as the transformer degradation component. 
DLMCs are reduced compared to the unconstrained case, but they still exhibit spikes compared to Full-opt.


\subsection{High EV-PV Penetration} \label{increased}

\begin{figure}[tb]
\centering
\includegraphics[width=3.4in]{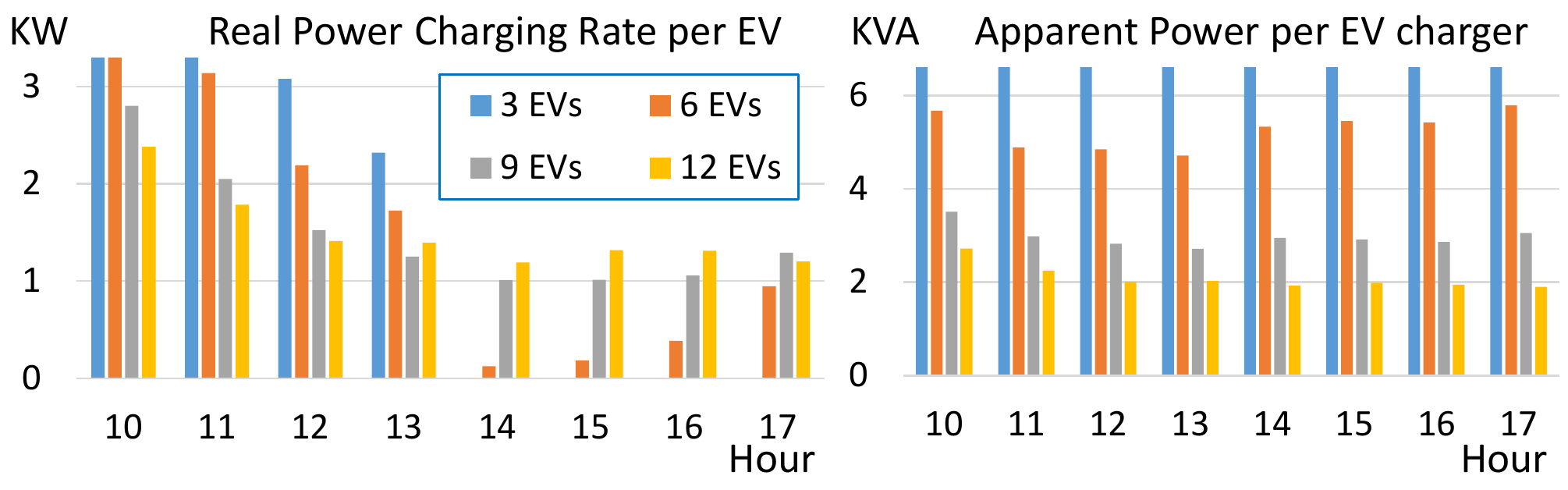}
\caption{EV charging rate (KW) and EV charger utilization (KVA) for 3, 6, 9, and 12-EV penetrations, under Full-opt, commercial node.}
\label{figEVhigh}
\end{figure}

We lastly tested scenarios with 9 and 12 EVs. 
Full-opt scheduling was still able to accommodate such high penetrations, whereas other scheduling options fail.
ToU keeps reasonably low aggregate LoL (40 hours) only for the 9-EV scenario at the residential node, although still 2-fold higher than Full-opt, whereas 32-fold higher LoL values are observed for the 12-EV penetration scenario.
At the commercial node, BaU/ToU exhibit aggregate LoL values that exceed 1000 hours (PQ-opt is even worse). 
Fig. \ref{figEVhigh} illustrates the charging rate and the utilization of the inverter that explains the robust performance of Full-opt with increasing EV penetrations.
Notably, for the 12-EV scenario, Full-opt reaches 114 aggregate LoL hours, which is still 40 times lower than the values encountered under BaU/ToU schedules.

Experimentation with 90 and 120 KVA PV penetration scenarios showed that Full-opt can significantly increase PV hosting capacity ``curtailing'' real and reactive power as needed, and, thus, avoiding excessive transformer temperatures and LoL.
Notably, if we consider self-scheduling PVs responding to optimal DLMC-based signals, Full-opt DLMCs support economically efficient ``curtailment.''


\section{Conclusions} \label{Conclusions}

We have shown the significance of distribution network asset degradation costs exemplified by transformer marginal LoL.
Our numerical results inquiry has focused on service transformers that are sensitive to granular/location-specific DER schedules.
A study focusing on the T\&D interface particularly during the few weeks following a station transformer failure, could and should model station transformers placing them on the line connecting the root node $0$ to the immediately downstream node $1$.
During periods of station transformer overloading, the LMP at node $0$ would be translated to a higher DLMC downstream of the station. 

We have shown that optimal DLMC-driven DER scheduling can achieve significantly higher hosting of DERs by distribution networks with minimal infrastructure investments.
Perhaps more importantly, we have shown that when DER schedules adapt to feeder DLMCs, the equilibrium DLMC profiles are rather smooth, thus allaying fears that DLMC volatility may lead to ratepayer revolt if DLMC-based rates are adopted. 
We proposed appropriate AC OPF models for granular marginal costing on distribution networks, and a real distribution network carried out on a feeder pilot study to obtain supporting numerical results. 
Comparison with popular open-loop DER scheduling options provided solid evidence that optimal DLMC-based clearing markets can bring about significant economic efficiencies and support the sustainability/adequacy of current distribution network infrastructure in the presence of high DER (PV, EV, and the like) adoption. 
It is particularly notable that the presence of DERs scheduled optimally through spatiotemporal DLMC adaptation will have a significant impact on grid asset/wires capacity expansion choices. Recent work (see \cite{LMV}) investigates the associated bridge between investment and operational planning.
\section*{Acknowledgment}
The authors would like to thank Holyoke Gas and Electric for providing actual feeder data for the pilot study.

\appendices

\section{Sensitivity Calculations} \label{AppB}

In what follows, we present the linear system associated with the calculations of the sensitivities required to unbundle the DLMC components.
The notation follows the AC OPF model introduced in \cite{PartI}.

Assume a system operating point that corresponds to network variables $\hat P_{ij,t}$, $\hat Q_{ij,t}$, $\hat v_{j,t}$, and $\hat l_{ij,t}$, $\forall j \in \mathcal{N}^+, t \in \mathcal{T}^+$.
We note that the partial derivatives of any network variable at time period $t$ \emph{w.r.t.} to net demand at time period $t'$ are all zero for $t \neq t'$, and that the partial derivatives of the voltage at the root node, $v_{0,t}$, are zero. 

Let us consider the partial derivatives \emph{w.r.t.} $p_{j',t'}$. 
At time period $t'$, the branch flow equations \cite[Eqs. (2)--(5)]{PartI} yield a $4N \times 4N$ system of linear equations:
\begin{equation} \label{SystemSens}
\left( \begin{matrix}
   {{\mathbf{A}}^{PP}} & {{\mathbf{0}}} & {{\mathbf{0}}} & {{\mathbf{A}}^{Pl}}  \\
   {{\mathbf{0}}} & {{\mathbf{A}}^{QQ}} & {{\mathbf{0}}} & {{\mathbf{A}}^{Ql}}  \\
   {{\mathbf{A}}^{vP}} & {{\mathbf{A}}^{vQ}} & {{\mathbf{A}}^{vv}} & {{\mathbf{A}}^{vl}}  \\
   {{\mathbf{A}}^{lP}} & {{\mathbf{A}}^{lQ}} & {{\mathbf{A}}^{lv}} & {{\mathbf{A}}^{ll}}  \\
\end{matrix} \right)\left( \begin{matrix}
   {{\mathbf{x}}^{P}}  \\
   {{\mathbf{x}}^{Q}}  \\
   {{\mathbf{x}}^{v}}  \\
   {{\mathbf{x}}^{l}}  \\
\end{matrix} \right)=\left( \begin{matrix}
   {{\mathbf{b}}^{P}}  \\
   {{\mathbf{b}}^{Q}}  \\
   {{\mathbf{0}}}  \\
   {{\mathbf{0}}}  \\
\end{matrix} \right), 
\end{equation}
where matrices $\mathbf{A}$ and vectors $\mathbf{x},\mathbf{b}$, are $N\times N$ matrices, and $N \times 1$ vectors, respectively, with $\mathbf{x}^P$,  $\mathbf{x}^Q$, $\mathbf{x}^v$, and $\mathbf{x}^l$ vectors of the partial derivatives $\frac{\partial P_{ij,t'}}{\partial p_{j',t'}} $,$\frac{\partial Q_{ij,t'}}{\partial p_{j',t'}}$, $\frac{\partial v_{j,t'}}{\partial p_{j',t'}}$, and $\frac{\partial l_{ij,t'}}{\partial p_{j',t'}}$, $\forall j \in \mathcal{N}^+$.
Let $\text{a}^{PP}_{jk}$ denote the $jk$-th element of matrix $\mathbf{A}^{PP}$, and $\text{b}_j$ the $j$-th element of vector $\mathbf{b}^P$. 
Similarly for the other matrices and vectors.
Also, let $i$ denote the node preceding $j$.
Then, we define matrices $\mathbf{A}$ and vectors $\mathbf{b}$ by providing the non-zero elements, $\forall j \in \mathcal{N}^+$, below.

For the first row representing \cite[Eq. (2)]{PartI}: $\text{a}^{PP}_{jj} = 1$,  $\text{a}^{PP}_{jk} = -1$, $\forall k: j \to k$, $\text{a}^{Pl}_{jj} = -r_{ij}$, and $\text{b}^{P}_{j'} = 1$. 
For the second row, representing \cite[Eq. (3)]{PartI}, $\text{a}^{QQ}_{jj} = 1$,  $\text{a}^{QQ}_{jk} = -1$, $\forall k: j \to k$, $\text{a}^{Ql}_{jj} = -x_{ij}$. 
For the third row, representing \cite[Eq. (4)]{PartI}, $\text{a}^{vP}_{jj} = 2 r_{ij}$,  $\text{a}^{vQ}_{jk} = 2 x_{ij}$, $\text{a}^{vv}_{ji} = -1$, (for $i \neq 0)$, $\text{a}^{vv}_{jj} = 1$, $\text{a}^{vl}_{jj} = - r_{ij}^2 - x_{ij}^2$. 
For the fourth row, representing \cite[Eq. (5)]{PartI}, $\text{a}^{lP}_{jj} = - 2 \hat P_{ij,t'}$, $\text{a}^{lQ}_{jj} = - 2 \hat Q_{ij,t'}$, $\text{a}^{lv}_{ji} = \hat l_{ij,t'}$ (for $i \neq 0)$, $\text{a}^{ll}_{jj} = \hat v_{i,t'}$.
Taking the partial derivatives \emph{w.r.t.} $q_{j',t'}$ in vectors $\mathbf{x}$, the only change is in the rhs of \eqref{SystemSens}.
Vector $\mathbf{b}^P$ should be zero, whereas $\mathbf{b}^Q$ has one non-zero element for $\text{b}^{Q}_{j'} = 1$.

\ifCLASSOPTIONcaptionsoff 
  \newpage
\fi


\begin{thebibliography}{99}

\bibitem{PartI}
P.~Andrianesis, and M.~Caramanis, 
``Distribution network marginal costs --- Part I: A novel AC OPF including transformer degradation,'' submitted, 2019.

\bibitem{CIGRE}
CIGRE, \emph{Benchmark systems for network integration of renewable and distributed energy resources}, Task Force C6.04, 2014.

\bibitem{BaiEtAl-DLMP}
L.~Bai, J.~Wang, C.~Wang, C.~Chen, and F.~Li, ``Distribution Locational Marginal Pricing (DLMP) for congestion management and voltage support,'' \emph{IEEE Trans. Power Syst.}, vol. 33, no. 4, pp. 4061--4073, 2018.

\bibitem{Abdelouadoud_EtAl_2015}
S.~Y.~Abdelouadoud, R.~Girard, F.~P~Neirac, and T.~Guiot, ``Optimal power flow of a distribution system based on increasingly tight cutting planes added to a second order cone relaxation,'' \emph{Int. J. Electr. Power Energy Syst.}, vol. 69, pp. 9--17, 2015.

\bibitem{HuangEtAl_2017}
S.~Huang, Q.~Wei, J.~Wang, and H.~Zhao, ``A sufficient condition on convex relaxation of AC optimal power flow in distribution networks,'' \emph{IEEE Trans. Power Syst.}, vol. 32, no. 2, pp. 1359--1368, 2017.

\bibitem{WeiEtAl_2017}
W.~Wei, J.~Wang, N.~Li, and S.~Mei, ``Optimal power flow of radial networks and its variations: A sequential convex optimization approach,'' \emph{IEEE Trans. Smart Grid}, vol. 8, no. 6, pp. 2974--2987, 2017.

\bibitem{NickEtAl_2018}
M.~Nick, R.~Cherkaoui, J.-Y.~Le~Boudec, and M.~Paolone, ``An exact convex formulation of the optimal power flow in radial distribution networks including transverse components,'' \emph{IEEE Trans. Autom. Control}, vol. 63, no. 3, pp. 682--697, 2018.

\bibitem{YuanEtAl-DLMP}
Z.~Yuan, M.~R.~Hesamzadeh, and D.~R.~Biggar, ``Distribution locational marginal pricing by convexified ACOPF and hierarchical dispatch,'' \emph{IEEE Trans. Smart Grid}, vol. 9, no. 4, pp. 3133--3142, 2018.

\bibitem{HilsheyEtAl_2013}
A.~D.~Hilshey, P.~D.~H.~Hines, P.~Rezaei, and J.~R.~Dowds, ``Estimating the impact of electric vehicle smart charging on distribution transformer aging,'' \emph{IEEE Trans. Smart Grid}, vol. 4, no. 2, pp. 905--913, 2013.
\bibitem{LMV}
P.~Andrianesis, M.~Caramanis, R.~Masiello, R.~Tabors, and S.~Bahramirad, ``Locational marginal value of distributed energy resources as non wires alternatives,'' \emph{IEEE Trans. Smart Grid}, accepted, 2019.







 










%

%
%
%

%
\end{thebibliography}
\end{document}